\documentclass[amsmath,10pt]{article}

\usepackage{graphics}
\usepackage[usenames]{color}
\usepackage{amsmath}
\usepackage{amsthm}
\usepackage{amsbsy}

\begin{document}

\bibliographystyle{plain}

\large \bf \begin{center}A UNIVERSAL CLOSED FORM FOR SQUARE MATRIX POWERS \normalsize

\rm Walter Shur
\end{center}

\rm \begin{center}20 Speyside Circle, Pittsboro, North Carolina, 27312, USA

 Tel: 919-542-7179 \ \ Email: wrshur@gmail.com

Affiliation: None (Retired) \end{center}\vspace{.2in}

\begin{center}\bf ABSTRACT.\end{center}

 \rm  This note presents a  simple, universal closed form for the powers of any square matrix. A diligent search of the internet gave no indication that the form is known.
\vspace{.4in}

\bf Keywords:\ \rm   Matrix, Powers, Closed Form\rm

\bf AMS Subject Classification:\ \rm 15A99, 65F30

\vspace{.3in}

1. \bf  THE CLOSED FORM \rm 

 Let $M=[m_{i,j}]$ be an n x n square matrix, real or complex, and let   $_nm_{i,j}$ be the $(i,j)^{th}$ element of $M^n$.  Then  $_nm_{i,j}$ can be expressed as the sum of $n$ terms, each of which is the product of a constant not depending on $n$, a binomial coefficient, and a power of an eigenvalue. In particular, for each unique eigenvalue, $\lambda$, of $M$, with multiplicity $mp_\lambda$, the portion of  $_nm_{i,j}$ arising from $\lambda$ is

$c_{i,j,1} \binom{n-1}{0}\lambda^{n-1}+c_{i,j,2}\binom{n-1}{1}\lambda^{n-2}+\cdots +c _{i,j,mp_\lambda}\binom{n-1}{mp_\lambda-1}\lambda^{n-mp_\lambda}$.

If the matrix is not singular, the formula works for any power of $M$, negative, zero or positive; if the matrix is  singular, for positive powers only. In the case of zero eigenvalues,  $0^n$ is taken as $0$ if $n<0$, and $1$ if $n=0$.

\pagebreak

\bf 2. ILLUSTRATION \rm

Let $M$ be the matrix\footnote{this matrix appeared in [3]}

 \[ M= \left(\begin{array}{cccccc}
66&83&95&31&-50&-63\\
-71&-79&-86&-22&59&72\\
83&74&69&4&-77&-90\\
-74&-50&-34&16&77&86\\
-7&-31&-47&-32&-12&-5\\
65&89&105&41&-40&-56
\end{array}\right)\]

The $6$ eigenvalues of $M$ are $\{-3,2,2,1,1,1\}$. Hence,

\hspace{-1.3in} \mbox{ $_nm_{i,j}=c_{i,j,1}(-3)^{n-1}+c_{i,j,2} 2^{n-1}+c_{i,j,3} \binom{n-1}{1}2^{n-2}+c_{i,j,4} 1^{n-1}+c_{i,j,5} \binom{n-1}{1}1^{n-2}+c_{i,j,6} \binom{n-1}{2}1^{n-3}$.}

To determine the values of the $c$ constants for say, $_nm_{2,5}$,  we first obtain by direct numerical calculation the values of $_nm_{2,5}$ for $n=1$  to  $6$, namely, $59, 229, 764, 1915, 5270$, and $11377$, respectively.

We then solve the $6$ simultaneous equations 

\hspace{-1.3in}\mbox{$_1m_{2,5}=59$, $_2m_{2,5}=229$, $_3m_{2,5}=764$, $_4m_{2,5}=1915$, $_5m_{2,5}=5270$, $_6m_{2,5}=11377$,}

 obtaining

\mbox{$c_{2,5,1}=3, c_{2,5,2}=65, c_{2,5,3}=126, c_{2,5,4}=-9, c_{2,5,5}=-9, c_{2,5,6}=0$.}

 Hence, we have the closed form

\hspace{-.7in} $_nm_{2,5}=3(-3)^{n-1}+65\  2^{n-1}+126\binom{n-1}{1}2^{n-2}-9 \ 1^{n-1}-9\binom{n-1}{1}1^{n-2}.$

\bf 3. DERIVATION OF THE CLOSED FORM \rm

Any square matrix $A$ can be written as $A = Q T Q^{-1}$,
where $Q$ is a unitary matrix, and $T$ is an upper triangular matrix (Schur decomposition). Hence, $A^n=Q T^n Q^{-1}$. The eigenvalues of $A$ and $T$ are the same, and are the main diagonal elements of $T$. Theorem 3 in [2] shows that every element of $T^n$ can be expressed by the linear combination described above in the first section of this note. Therefore, every element of $A^n$ can also be expressed by that same general form.

\bf 4. CONCLUSION \rm

The purpose of this note was not to present a new method of calculating the numerical values of the elements of a matrix raised to the $n^{th}$ power. Nor does it present the only way to obtain a closed form for matrix powers. For example, the Mathematca\footnote{Version 10.3.0.0} function MatrixPower[m,n][[2,5]] gives

$_nm_{2,5}=\frac{1}{2}(-2(-3)^n+2^{1+n}-18n+63\  2^n n)$. 

This is identical to the closed form for $_nm_{2,5}$ given above, but with the beauty and symmetry absent.

Rather, the primary purpose of this note was to share with others the aesthetic beauty of a simple universal form, which represents some very complicated calculations.

\vspace{.1in}
\hspace{4in} 11/30/15

\end{document}